\documentclass[12pt, draft]{article}
\usepackage{amsmath, amsfonts, amsthm, amssymb}
\newtheorem{theorem}{Theorem}[section]

\newtheorem{Proposition}[theorem]{Proposition}
\newtheorem{corollary}[theorem]{Corollary}
\newtheorem{definition}{Definition}[subsection]

\begin{document}
\title{\bf THE TOPOLOGY OF FOLIATIONS FORMED BY
THE GENERIC K-ORBITS OF A SUBCLASS OF THE INDECOMPOSABLE MD5-GROUPS}
\author{${\bf Le Anh Vu}^{\dag}$ \quad and \quad {\bf Duong Quang Hoa}\\
$\dag$: Department of Mathematics and Informatics\\
University of Pedagogy, Ho Chi Minh City, Vietnam\\E-mail:\,
vula@math.hcmup.edu.vn} \footnotetext{

2000AMS Mathematics Subject Classification: Primary 22E45, Secondary
46E25, 20C20.}
\date{}
\maketitle

\begin{abstract}
    The present paper is a continuation of [13], [14] of the authors. Specifically,
    the paper considers the MD5-foliations associated to connected and
    simply connected MD5-groups such that their Lie algebras have
    4-dimensional commutative derived ideal. In the paper, we give the
    topological classification of all considered MD5-foliations. A description
    of these foliations by certain fibrations or suitable actions of $\mathbb{R}^{2}$
    and the Connes' C*-algebras of the foliations which come from
    fibrations are also given in the paper.
\end{abstract}

\subsection*{Introduction}
Our point of departure is the problem of finding out of the classes
of C*-algebras which can be described by the operator KK-functors.
In 1980, studying the Kirillov's Orbit Method, D. N. Diep suggested
to consider the class of MD-groups. An MD-group with dimension n
(for brevity, an MDn-group) in his terms (see [2, Section 4.1]) is
an n-dimensional solvable real Lie group whose orbits in the
co-adjoint representation (i.e. the K-representation) are the orbits
of zero or maximal dimension. The Lie algebra of each MDn-group is
called an MDn-algebra.

The first reason for studying the MD-groups is the fact that the
C*-algebras of the MD-groups can be described by KK-functors (see
[2, Chapters 3, 5], [4, Section 3]). On the other hand, for every
MD-group G, the family of K-orbits of maximal dimension forms a
measured foliation in terms of A. Connes ([1]). This foliation is
called MD-foliation associated with G. In general, the leaf space of
a foliation with the quotient topology is a fairly untractable
topological space. To improve upon the shortcoming, A. Connes
associates with each measurable foliation a C*-algebra (see [1,
Section 2]). In the cases of Reeb foliations, the method of
KK-functor has been proved very effective in describing the Connes'
C*-algebras by A. M. Torpe (see [5, Sections 3, 4]).

Combining methods of A. Kirillov (see [3, Section 15]) and A. Connes
(see [1, Sections 2, 5]),  Vu has considered MD4-foliations
associated with all indecomposable connected MD4-groups in [6], [7],
[8]. Recently, in [9], [11], [12] Vu together with Nguyen Cong Tri
anh Duong Minh Thanh have studied  MD5-foliations associated with
the MD5-groups such that the first derived ideal of corresponding
Lie algebras is $\mathbb{R}^{k};\, k = 1, 2, 3$. The present paper
is a continuation of [13], [14] of the authors and it is concerned
with MD5-foliations associated with the indecomposable connected and
simply connected MD5-groups such that MD5-algebras of them have
4-dimensional commutative derived ideal.

We shall begin our discussion in Section 1 by recalling some
preliminary results and notations which will be used later. For more
details we refer the reader to References [1], [3].

In Section 2 we shall relist all MD5-algebras with 4-dimensional
commutative derived ideal and recall the geometric description of
the K-orbits of corresponding MD5-groups which have been announced
in [13], [14].

Section 3 will be devoted to the discussion of the MD5-foliations
associated with considered MD5-groups. We shall give a topological
classification of these foliations and describe them by certain
fibrations or suitable actions of $\mathbb{R}^{2}$. In addition, the
Connes' C*-algebras of the MD5-foliations which come from certain
fibrations will be also analytically characterized.

\section{Preliminaries}

\subsection{Foliations}

Let V be a smooth manifold. We denote by TV its tangent bundle, so
that for each x $\in$ V, $T_{x}V$ is the tangent space of V at x.

\begin{definition} A smooth subbundle $\mathcal{F}$ of TV is called
integrable if and only if the following condition is satisfied:
every $x \in V$ is contained in a submanifold W of V such that
$T_{p}$(W) = ${\mathcal{F}}_{p}\, (\forall p \in
W)$.\end{definition}

\begin{definition} A foliation (V, $\, \mathcal {F}$ ) is
given by a smooth manifold V and an integrable subbundle
$\mathcal{F}$ of TV. Then, V is called the foliated manifold and
$\mathcal{F}$ is called the subbundle defining the foliation. The
dimension of $\mathcal{F}$ is also called the dimension of foliation
(V, $\mathcal{F}$).\end{definition}

\begin{definition} Each maximal connected submanifolds L of V such
that $T_{x}(L) = {\mathcal{F}}_{x}$\,($\forall x \in L$) is called a
leaf of the foliation (V, $\, \mathcal{F}$).\end{definition}

    The set of leaves with the quotient topology is denoted by
    V/${\mathcal{F}}$ and called the \emph{space of leaves} or \emph{leaf space} of
    (V, $\, \mathcal{F}$). In general, it is a fairly untractable topological space.

   The partition of V in leaves: V = $\bigcup_{\alpha \in V/\mathcal{F}}L_{\alpha}$
   is characterized geometrically by the following local triviality: Every
   $x \in V$ has a system of local coordinates  $\{ U; x^{1}, x^{2}, ..., x^{n} \}
   (x \in U; n = dim V)$  so that for any leaf
   $L$ with $L \cap$ U $\ne \emptyset$ , each connected component of
   $L \cap U$ (which is called a \emph{plaque} of the leaf $L$)
   is given by the equations
$$x^{k+1} = c^{1},\, x^{k+2} = c^{2},\, ...\,, x^{n}\, =\, c^{n-k};\, k = dim \mathcal{F}< n;$$
where $c^{1}, c^{2}, ..., c^{n-k}$ are constants (depending on each
plaque). Each such system $\{ U, x^{1}, x^{2}, ..., x^{n} \}$ is
called a \emph{foliation chart}.

   A foliation can be given by a partition of V in a family
$\mathcal{C}$ of its submanifolds if there exist some integrable
subbundle $\mathcal{F}$ of TV such that each $L\in{\mathcal{C}}$ is
a maximal connected integral submanifold of . Then $\mathcal{C}$ is
the family of leaves of the foliation (V, $\, \mathcal{F}$).
Sometimes $\mathcal{C}$ is identified with $\mathcal{F}$ and we say
that (V, $\, \mathcal{F}$) is formed by $\mathcal{C}$.

\subsection{Measurable Foliations}

\begin{definition} A submanifold N of the foliated manifold V is
called a transversal if and only if $T_{x}V = T_{x}N \oplus
{\mathcal{F}}_{x},\, \forall x \in N$. Thus,
 dimN = n - dim$\mathcal{F}$ = codim$\mathcal{F}$.

    A Borel subset B of V such that $B\cap L$ is countable for any leaf L
    is called a Borel transversal to
    (V,\,$\mathcal{F}$).\end{definition}

\begin{definition} A transverse measure $\Lambda$ for the
foliation (V,\,$\mathcal{F}$) is $\sigma$ - additive map B\,
$\mapsto \Lambda$ (B) from the set of all Borel transversals to [0,
+$\infty$] such that the following conditions are satisfied :

    (i) If $\psi$ : $B_{1} \rightarrow B_{2}$ is a Borel bijection and
    $\psi$(x) is on the leaf of any x$\in B_{1}$, then $\Lambda(B_{1}) = \Lambda(B_{2})$.

    (ii) $\Lambda(K)< + \infty$ if K is any compact subset of a smooth transversal
    submanifold of V.

    By a measurable foliation we mean a foliation
    (V,\,$\mathcal{F}$) equipped with some transverse measure
    $\Lambda$.\end{definition}

 Let (V,\,$\mathcal{F}$) be a foliation with $\mathcal{F}$ is
oriented. Then the complement of zero section of the bundle
${\Lambda}^{k}(\mathcal{F})$ (k = dim$\mathcal{F}$)  has two
components ${\Lambda}^{k}{(\mathcal{F})}^{-}$ and
${\Lambda}^{k}{(\mathcal{F})}^{+}$.

    Let $\mu$ be a measure on V and  $\{ U, x^{1}, x^{2}, ..., x^{n} \}$
    be a foliation chart of (\nobreak V,\,$\mathcal{F}$\nobreak).
    Then U can be identified with the direct product $N \times {\Pi}$
    of some smooth transversal submanifold N of V and some plaque
$\Pi$. The restriction of $\mu$ on $U \equiv N \times {\Pi}$ becomes
the product ${\mu}_{N} \times {\mu}_{\Pi}$ of measures ${\mu}_{N}$
and ${\mu}_{\Pi}$ respectively.

     Let X $\in C^{\infty}{\bigl({\Lambda}^{k}(\mathcal{F})\bigr)}^{+}$
     be a smooth k-vector field and ${\mu}_{X}$ be the measure on each leaf L
     determined by the volume element X.

\begin{definition} The measure $\mu$ is called
X-invariant if and only if ${\mu}_{X}$ is proportional to
${\mu}_{\Pi}$ for an arbitrary foliation chart $\{ U, x^{1}, x^{2},
..., x^{n} \}$.\end{definition}

 Let (X, $\mu$), (Y, $\nu$) be two pairs where X,Y
 $\in C^{\infty}{\bigl({\Lambda}^{k}(\mathcal{F})\bigr)}^{+}$ and $\mu, \nu$
 are measures on V such that $\mu$ is X-invariant, $\nu$ is Y-invariant.

\begin{definition} ( X, $\mu$ ), ( Y, $\nu$ ) are
equivalent if and only if Y = $\varphi$ X and $\mu$ = $\varphi \nu$
for some $\varphi \in C^{\infty}(V).$\end{definition}

    There is one bijective map between the set of transverse measures for
(V,\, $\mathcal{F}$) and the set of equivalent classes of pairs (X,
\, $\mu$), where X $\in
C^{\infty}{\bigl({\Lambda}^{k}(\mathcal{F})\bigr)}^{+}$ and $\mu$ is
a X-invariant measure on V.

    Thus, to prove that (V,\,$\mathcal{F}$) is measurable, we only need to choose
    some suitable pair (X, $\mu$) on V.

\section{MD5-algebras with 4-dimensional Commutative Derived Ideal
and Geometry of K-orbits of Corresponding MD5-Groups}

\subsection{The list of MD5-algebras with 4-dimensional Commutative Derived Ideal}
 In [13] the first author have listed all MD5-algebras such that
 their first derived ideals are $\mathbb{R}^{4}$. For completeness,
 we will relist all these MD5-algebras here.

\begin{Proposition}[{see [13, Theorem 3.2]}] Let $\mathcal{G}$ be an indecomposable
MD5-algebra with $\mathcal{G}^{1}:= [\mathcal{G}, \mathcal{G}]\cong
\mathbb{R}^{4}$. Then we can choose a suitable basis $( X_{1},
X_{2},\\ X_{3}, X_{4}, X_{5} )$ of \, $\mathcal{G}$ such that
${\mathcal{G}}^{1} = \mathbb{R}.X_{3} \oplus \mathbb{R}.X_{3} \oplus
\mathbb{R}.X_{4} \oplus \mathbb{R}.X_{5} \equiv {\mathbb{R}}^{4}$,\,
$ad_{X_{1}} \in End({\mathcal{G}}^{1}) \equiv Mat_{4}(\mathbb{R})$
and $\mathcal{G}$ is isomorphic to one and only one of the following
Lie algebras.
            \begin{description}
                \item[1.]${\mathcal{G}}_{5,4,1({\lambda}_{1}, {\lambda}_{2}, {\lambda}_{3})}:$
                    $$ad_{{X}_1} = \begin{pmatrix} {{\lambda}_1}&0&0&0\\
                    0&{{\lambda}_2}&0&0\\0&0&{\lambda}_{3}&0\\0&0&0&1\end{pmatrix};$$
                    $${\lambda}_1, {\lambda}_2, {\lambda}_3 \in \mathbb{R}\setminus
                    \lbrace 0, 1\rbrace,\quad
                    {\lambda}_1 \neq {\lambda}_2 \neq {\lambda}_3 \neq {\lambda}_1.$$ \vskip 0.5cm
                \item[2.]${\mathcal{G}}_{5,4,2({\lambda}_{1}, {\lambda}_{2})}:$
                    $$ad_{{X}_1} = \begin{pmatrix} {{\lambda}_1}&0&0&0\\
                    0&{{\lambda}_2}&0&0\\0&0&1&0\\0&0&0&1\end{pmatrix};\quad
                    {\lambda}_{1}, {\lambda}_{2} \in \mathbb{R}\setminus \lbrace 0, 1 \rbrace ,
                    {\lambda}_1 \neq {\lambda}_2.$$ \vskip 0.5cm
                \item[3.]${\mathcal{G}}_{5,4,3(\lambda)}:$
                    $$ad_{{X}_1} = \begin{pmatrix}
                    {\lambda}&0&0&0\\0&{\lambda}&0&0\\0&0&1&0\\0&0&0&1 \end{pmatrix}; \quad
                    {\lambda} \in \mathbb{R}\setminus \lbrace 0, 1 \rbrace .$$ \vskip 0.5cm
                \item[4.]${\mathcal{G}}_{5,4,4(\lambda)}:$
                    $$ad_{{X}_1} = \begin{pmatrix} {\lambda}&0&0&0\\0&1&0&0\\
                    0&0&1&0\\0&0&0&1 \end{pmatrix};\quad
                    {\lambda} \in \mathbb{R}\setminus \lbrace 0, 1 \rbrace.$$ \vskip 0.5cm
                \item[5.]${\mathcal{G}}_{5,4,5}:$
                    $$ad_{{X}_1} = \begin{pmatrix} 1&0&0&0\\0&1&0&0\\
                    0&0&1&0\\0&0&0&1 \end{pmatrix}.$$ \vskip 0.5cm
                \item[6.]${\mathcal{G}}_{5,4,6({\lambda}_{1}, {\lambda}_{2})}$ :
                    $$ad_{{X}_1} = \begin{pmatrix} {{\lambda}_1}&0&0&0\\
                    0&{{\lambda}_2}&0&0\\0&0&1&1\\0&0&0&1\end{pmatrix};\quad
                    {\lambda}_{1}, {\lambda}_{2} \in \mathbb{R}\setminus \lbrace 0, 1 \rbrace ,
                    {\lambda}_1 \neq {\lambda}_2.$$ \vskip 0.5cm
                \item[7.]${\mathcal{G}}_{5,4,7(\lambda)}:$
                    $$ad_{{X}_1} = \begin{pmatrix}
                    {\lambda}&0&0&0\\0&{\lambda}&0&0\\0&0&1&1\\0&0&0&1 \end{pmatrix}; \quad
                    {\lambda} \in \mathbb{R}\setminus \lbrace 0, 1 \rbrace .$$ \vskip 0.5cm
                \item[8.]${\mathcal{G}}_{5,4,8(\lambda)}:$
                    $$ad_{{X}_1} = \begin{pmatrix}
                    {\lambda}&1&0&0\\0&{\lambda}&0&0\\0&0&1&1\\0&0&0&1 \end{pmatrix}; \quad
                    {\lambda} \in \mathbb{R}\setminus \lbrace 0, 1 \rbrace .$$ \vskip 0.5cm
                \item[9.]${\mathcal{G}}_{5,4,9(\lambda)}:$
                    $$ad_{{X}_1} = \begin{pmatrix}
                    {\lambda}&0&0&0\\0&1&1&0\\0&0&1&1\\0&0&0&1 \end{pmatrix}; \quad
                    {\lambda} \in \mathbb{R}\setminus \lbrace 0, 1 \rbrace .$$ \vskip 0.5cm
                \item[10.]${\mathcal{G}}_{5,4,10}:$
                    $$ad_{{X}_1} = \begin{pmatrix} 1&1&0&0\\0&1&1&0\\
                    0&0&1&1\\0&0&0&1 \end{pmatrix}.$$ \vskip 0.5cm
                \item[11.]${\mathcal{G}}_{5,4,11({\lambda}_{1}, {\lambda}_{2},\varphi)}:$
                    $$ad_{{X}_1} = \begin{pmatrix} cos\varphi&-sin\varphi&0&0\\
                    sin\varphi&cos\varphi&0&0\\0&0&{\lambda}_{1}&0\\0&0&0&{\lambda}_{2}\end{pmatrix};$$
                    $${\lambda}_{1}, {\lambda}_{2} \in \mathbb{R}\setminus \lbrace 0 \rbrace ,
                    {\lambda}_1 \neq {\lambda}_2,\varphi \in (0,\pi).$$ \vskip 0.5cm
                \item[12.]${\mathcal{G}}_{5,4,12(\lambda, \varphi)}:$
                    $$ad_{{X}_1} = \begin{pmatrix} cos\varphi&-sin\varphi&0&0\\
                    sin\varphi&cos\varphi&0&0\\0&0&\lambda&0\\0&0&0&\lambda\end{pmatrix};\quad
                    \lambda \in \mathbb{R}\setminus \lbrace 0 \rbrace, \varphi \in (0,\pi).$$ \vskip 0.5cm
                \item[13.]${\mathcal{G}}_{5,4,13(\lambda, \varphi)}:$
                    $$ad_{{X}_1} = \begin{pmatrix} cos\varphi&-sin\varphi&0&0\\
                    sin\varphi&cos\varphi&0&0\\0&0&\lambda&1\\0&0&0&\lambda\end{pmatrix};\quad
                    \lambda \in \mathbb{R}\setminus \lbrace 0 \rbrace, \varphi \in (0,\pi).$$ \vskip 0.5cm
                \item[14.]${\mathcal{G}}_{5,4,14(\lambda, \mu, \varphi)}:$
                    $$ad_{{X}_1} = \begin{pmatrix} cos\varphi&-sin\varphi&0&0\\
                    sin\varphi&cos\varphi&0&0\\0&0&\lambda&-\mu\\0&0&\mu&\lambda\end{pmatrix};$$
                    $$\lambda, \mu \in \mathbb{R}, \mu > 0, \varphi \in
                    (0,\pi).$$ \hfill{$\square$}
            \end{description}
\end{Proposition}

\subsection*{Remarks} Let us recall that each real Lie
algebra $\mathcal{G}$ define only one connected and simply connected
Lie group G such that Lie(G) = $\mathcal{G}$. Therefore we obtain a
collection of fourteen families of connected and simply connected
MD5-groups corresponding to the indecomposable MD5-algebras given in
Theorem 2.1. For convenience, each MD5-group from this collection is
also denoted by the same indices as corresponding MD5-algebra. For
example, $G_{5,4,2({\lambda}_{1}, {\lambda}_{2})}$ is the connected
and simply connected MD5-group corresponding to
${\mathcal{G}}_{5,4,2({\lambda}_{1}, {\lambda}_{2})}$. All of these
groups are indecomposable MD5-groups. In [14], we have described
geometry of the K-orbits of them. But we now recall this result in
the next subsection for completeness.

\subsection{The Picture of K-orbits of Corresponding Connected and Simply Connected Groups}

Let G be one of considered MD5-groups, ${{\mathcal{G}} = < X_{1},
X_{2}, X_{3}, X_{4}, X_{5}>}$ be its Lie algebra,
${\mathcal{G}}^{*}$ = $<X_{1}^{*}$, $X_{2}^{*}$, $X_{3}^{*}$,
$X_{4}^{*}$, $X_{5}^{*}> \equiv \mathbb{R}^{5}$ be the dual space of
${\mathcal{G}}$, $F = {\alpha}X_{1}^{*} + {\beta}X_{2}^{*} +
{\gamma}X_{3}^{*} + {\delta}X_{4}^{*} + {\sigma}X_{5}^{*} \equiv
({\alpha}, {\beta}, {\gamma}, {\delta}, {\sigma})$ be an arbitrary
element of ${\mathcal{G}}^{*}$. The notation ${\Omega}_{F}$ will be
used to denote the K-orbit of G which contains $F$. The geometrical
picture of the K-orbits of considered MD5-groups is given by the
following proposition which has proved in [14]. \vskip0.8cm
\begin{Proposition}[{see[14, Theorems 3.3.1, 3.3.2, 3.3.3, 3.3.4]}] The
K-orbit ${\Omega}_{F}$ of G is described as follows.
    \begin{description}
        \item [1.] Let $G$ is one of \, $G_{5,4,1({\lambda}_{1}, {\lambda}_{2},
{\lambda}_{3})}$,\, $G_{5,4,2({\lambda}_{1}, {\lambda}_{2})}$,\,
$G_{5,4,3(\lambda)}$,\, $G_{5,4,4(\lambda)}$,\, $G_{5,4,5)}$,\\
$G_{5,4,6({\lambda}_{1}, {\lambda}_{2})}$, $G_{5,4,7(\lambda)}$,
$G_{5,4,8(\lambda)}$, $G_{5,4,9(\lambda)}$, $G_{5,4,10}$;
${\lambda}_{1}$, ${\lambda}_{2}$, ${\lambda}_{3},\lambda\in
\mathbb{R}\backslash\{0,1\}$.
        \begin{description}
            \item[1.1.] If $\beta=\gamma=\delta=\sigma=0$ then $\Omega_F=\{F\}$ (the
0-dimensional orbit).
            \item[1.2.] If $\beta^2+\gamma^2+\delta^2+\sigma^2\neq0$ then $\Omega_F$
is the orbit of dimension 2 and it is one of the following:
        \end{description}
\begin{itemize}
    \item $\{(x,\beta{e^{a\lambda_1}},\gamma{e^{a\lambda_2}},\delta{e^{a\lambda_3}},
    \sigma{e^a}),x,a\in\mathbb{R}\}$ when $G=G_{5,4,1({\lambda}_{1},
    {\lambda}_{2}, {\lambda}_{3})}$.
    \item $\{(x,\beta{e^{a\lambda_1}},\gamma{e^{a\lambda_2}},\delta{e^{a}},
    \sigma{e^a}),x,a\in\mathbb{R}\}$ when $G=G_{5,4,2({\lambda}_{1},
    {\lambda}_{2})}$.
    \item $\{(x,\beta{e^{a\lambda}},\gamma{e^{a\lambda}},\delta{e^a},
    \sigma{e^a}),x,a\in\mathbb{R}\}$ when $G=G_{5,4,3(\lambda)}$.
    \item $\{(x,\beta{e^{a\lambda}},\gamma{e^a},\delta{e^{a}},
    \sigma{e^a}),x,a\in\mathbb{R}\}$ when $G=G_{5,4,4(\lambda)}$.
    \item $\{(x,\beta{e^a},\gamma{e^a},\delta{e^a},
    \sigma{e^a}),x,a\in\mathbb{R}\}$ when $G=G_{5,4,5}$.
    \item $\{(x,\beta{e^{a\lambda_1}},\gamma{e^{a\lambda_2}},\delta{e^{a}},\delta{a}{e^{a}}+
    \sigma{e^a}),x,a\in\mathbb{R}\}$ when $G=G_{5,4,6({\lambda}_{1},{\lambda}_{2})}$.
    \item $\{(x,\beta{e^{a\lambda}},\gamma{e^{a\lambda}},\delta{e^a},\delta{a}{e^{a}}+
    \sigma{e^a}),x,a\in\mathbb{R}\}$ when $G=G_{5,4,7(\lambda)}$.
    \item $\{(x,\beta{e^{a\lambda}},\beta{a}{e^{a\lambda}}+\gamma{e^{a\lambda}},
    \delta{e^a},\delta{a}{e^{a}}+\sigma{e^a}),x,a\in\mathbb{R}\}$

    \hfill when $G=G_{5,4,8(\lambda)}$.
    \item $\{(x,\beta{e^{a\lambda}},\gamma{e^a},\gamma{a}{e^a}+\delta{e^{a}},
    \frac{\gamma{a^2}e^a}{2}+\delta{a}{e^{a}}+\sigma{e^a}),x,a\in\mathbb{R}\}$

    \hfill {when $G=G_{5,4,9(\lambda)}$}.
    \item $\{(x,\beta{e^a},\beta{a}{e^a}+\gamma{e^a},\frac{\beta{a^2}e^a}{2}+
    \gamma{a}{e^a}+\delta{e^{a}},\frac{\beta{a^3e^a}}{6}+ \frac{\gamma{a^2}e^a}{2}+
    \delta{a}{e^{a}}+\sigma{e^a})$,

    \hfill {$x,a\in\mathbb{R}\}$} when $G=G_{5,4,10}$.
\end{itemize}
    \item [2.] Let G is one of $G_{5,4,11({\lambda}_{1}, {\lambda}_{2},\varphi)}$,
$G_{5,4,12(\lambda,\varphi)}$,
$G_{5,4,13(\lambda,\varphi)};\lambda_1,\lambda_2,\lambda\in\mathbb{R}\setminus\{0\};\\
\varphi\in(0,\pi)$. Let us identify
${\mathcal{G}}_{5,4,11(\lambda_{1}, \lambda_{2},\varphi)}^*$,
${\mathcal{G}}_{5,4,12(\lambda,\varphi)}^*$,
${\mathcal{G}}_{5,4,13(\lambda,\varphi)}^*$ with $\mathbb{R}\times
\mathbb{C}\times {\mathbb{R}}^2$ and $F$ with
$(\alpha,\beta+i\gamma,\delta,\sigma)$. Then we have
    \begin{description}
    \item[2.1.] If
$\beta+i\gamma=\delta=\sigma=0$ then $\Omega_F=\{F\}$ (the
0-dimensional orbit).
    \item[2.2.] If $|\beta+i\gamma|^2+\delta^2+\sigma\neq0$ then $\Omega_F$
is the orbit of dimension 2 and it is one of the following:
    \end{description}
\begin{itemize}
    \item $\{(x,(\beta+i\gamma)e^{ae^{-i\varphi}},\delta{e^{a\lambda_1}},
    \sigma{e^{a\lambda_2}}),x,a\in\mathbb{R}\}$
    when $G=G_{5,4,11({\lambda}_{1}, {\lambda}_{2},
\varphi)}$.
    \item $\{(x,(\beta+i\gamma)e^{ae^{-i\varphi}},\delta{e^{a\lambda}},
    \sigma{e^{a\lambda}}),x,a\in\mathbb{R}\}$
    when $G=G_{5,4,12(\lambda,\varphi)}$.
    \item $\{(x,(\beta+i\gamma)e^{ae^{-i\varphi}},\delta{e^{a\lambda}},
    \delta{a}{e^{a\lambda}}+\sigma{e^{a\lambda}}),x,a\in\mathbb{R}\}$ when \\
    $G=G_{5,4,13(\lambda,\varphi)}$.
\end{itemize}
    \item [3.] Let $G$ is
$G_{5,4,14(\lambda, \mu, \varphi)}$. Let us identify
${\mathcal{G}}_{5,4,14(\lambda,\mu,\varphi)}^*$ with
$\mathbb{R}\times \mathbb{C} \times \mathbb{C}$ and $F$ with
$(\alpha,\beta+i\gamma,\delta+i\sigma)$;
$\lambda,\mu\in\mathbb{R};\mu>0;\varphi\in(o,\pi)$. Then we have
\begin{description}
    \item[3.1.] If $\beta+i\gamma=\delta+i\sigma=0$ then $\Omega_F=\{F\}$ (the
0-dimensional orbit).
    \item[3.2.] If $|\beta^2+i\gamma^2|+|\delta^2+i\sigma^2|\neq0$ then

$$\Omega_F=\{(x,(\beta+i\gamma)e^{ae^{-i\varphi}},(\delta+i\sigma)e^{a(\lambda-i\mu)}),x,a\in\mathbb{R}\}$$
(the 2-dimensional orbit). \hfill{$\square$}
\end{description}
\end{description} \end{Proposition}

\section{On MD5-foliations Associated
to\\ Considered MD5-groups and Their \\ Topological Classification}

\subsection{Foliations formed by K-orbits of dimension two of
considered MD5-groups}

In the introductory section we have emphasized that the family of
maximal-dimensional K-orbits of every MD-group forms a measured
foliation in terms of A. Connes. Namely, this fact has proved in
[7], [9], [11], [12]  not only  for all connected MD4-groups but
also for all connected MD5-groups with corresponding Lie algebras
have commutative derived ideal of dimension $k$ with $k < 4$. The
following is a similar assertion for the MD5-groups considered in
this paper. It is also proved by the same method as the one in [7],
[9], [11], [12]. So we omit the proof.

\begin{theorem} Let G be one of the connected and simply connected
MD5-groups corresponding to the MD5-algebras listed in Theorem 2.1,
$\mathcal{F}_{G}$ be the family of all its K-orbits of dimension two
and $V_{G}: = \bigcup \{ \Omega / \Omega \in \mathcal{F}_{G}\}$.
Then $(V_{G},\,\mathcal{F}_{G})$ is a measurable foliation in the
sense of Connes. We call it MD5-foliation associated with MD5-group
G. \hfill {$\square$}
\end{theorem}

\subsection*{Remarks and Notations}

It should be noted that $V_{G}$ is an open submanifold of the dual
space $\mathcal{G}^{*} \cong \mathbb{R}^{5}$ of the Lie algebra
$\mathcal{G}$ corresponding to G. Furthermore, for all MD5-groups of
the forms $G_{5,4, ...}$, the manifolds $V_{G}$ are diffeomorphic to
each other. So, for simplicity of notation, we shall write $(V,\,
F_{4, ...})$ instead of $(V_{G_{4, ...}},\, F_{G_{4, ...}})$. The
following theorem will be fundamental in this paper.

\begin{theorem}(The topological classification of considered
MD5-foliations)
\begin{enumerate}
    \item [1.]There exist exactly 3 topological types of fourteen
    families of considered MD5-foliations as follows:
    \begin{enumerate}
        \item[1.1.]
        $\Bigl \{(V,{\mathcal{F}}_{4,1(\lambda_1,\lambda_2,\lambda_3)}),
        (V,{\mathcal{F}}_{4,2(\lambda_1,\lambda_2)}),
        (V,{\mathcal{F}}_{4,3(\lambda)}), (V,{\mathcal{F}}_{4,4(\lambda)}),
        (V,{\mathcal{F}}_{4,5})$,\\
        $(V,{\mathcal{F}}_{4,6(\lambda_1,\lambda_2)}),
        (V,{\mathcal{F}}_{4,7(\lambda)}), (V,{\mathcal{F}}_{4,8(\lambda)}),
        (V,{\mathcal{F}}_{4,9(\lambda)}), (V,{\mathcal{F}}_{4,10})$;
        $\lambda,\lambda_1,$

         \hfill $\lambda_2, \lambda_3 \in \mathbb{R} \backslash\{0,1\}$ $\Bigl \}$.
        \item[1.2.]
        $\Bigl
        \{(V,{\mathcal{F}}_{4,11(\lambda_1,\lambda_2,\varphi)}),\,
        (V,{\mathcal{F}}_{4,12(\lambda,\varphi)}),\,
        (V,{\mathcal{F}}_{4,13(\lambda,\varphi)})$;\,

        \hfill $\lambda,\,\,\, \lambda_1,
        \lambda_2 \in \mathbb{R}\backslash\{0,1\}$; $\varphi\in(0,\pi)\Bigl \}$.
        \item[1.3.]
        $\Bigl \{(V,{\mathcal{F}}_{4,14(\lambda,\mu,\varphi)});\quad
        \mu,\lambda\in \mathbb{R},\, \mu>0,\, \varphi\in(0,\pi)\Bigl \}$.
    \end{enumerate}
    We denote these types by \, ${\mathcal{F}}_1,\, {\mathcal{F}}_2,\,
    {\mathcal{F}}_3$ respectively.
    \item [2.] Furthermore, we have
    \begin{enumerate}
        \item[2.1.] The MD5-foliations of types ${\mathcal{F}}_1$ are
        trivial fibration with connected fibres on the 3-dimensional
        unitary sphere $\mathbf{S}^3$.
        \item[2.2.] The MD5-foliations of types ${\mathcal{F}}_2,\,
        {\mathcal{F}}_3$ can be given by suitable actions of
        $\mathbb{R}^{2}$ on the foliated manifolds
        $V \cong{(\mathbb{R}^4)}^{*} \times \mathbb{R}$.
    \end{enumerate}
\end{enumerate}
\end{theorem}
\subsection*{Prove of Theorem 3.2}
Let us recall that two foliations $(V,{\mathcal{F}}),
(V,{\mathcal{F'}})$ are said to be topologically equivalent if there
exists a homeomorphism $h:V\rightarrow{V}$ which takes leaves of
$\mathcal{F}$ onto leaves of $\mathcal{F'}$. The map $h$ is called a
topological equivalence of considered foliations.
\begin{description}
    \item[1.] Firstly, we prove Assertion 1 in the
    theorem. Namely, we need to give the topological classification
    of considered MD5-foliations.
        \begin{enumerate}
        \item[1.1.] We consider maps
        $h_{4,1(\lambda_1,\lambda_2,\lambda_3)},
        h_{4,2(\lambda_1,\lambda_2)}, h_{4,3(\lambda)},
        h_{4,4(\lambda)}, h_{4,6(\lambda_1,\lambda_2)},\\
        h_{4,7(\lambda)}, h_{4,8(\lambda)}, h_{4,9(\lambda)},
        h_{4,10(\lambda)}$ from $V\cong \mathbb{R}\times{({\mathbb{R}}^4)}^{*}$ to $V$
        which are defined as follows.
\vskip0.7cm
        $h_{4,1(\lambda_1,\lambda_2,\lambda_3)}(x,y,z,t,s)$: =\,

                \hskip2cm $ = \,(x,sign(y).|y|^{\frac{1}{\lambda_1}},
                sign(z).|z|^{\frac{1}{\lambda_2}},
                sign(t).|t|^{\frac{1}{\lambda_3}},s).$
\vskip0.5cm
        $h_{4,2(\lambda_1,\lambda_2)}(x,y,z,t,s)$: =\,

                \hskip2cm $=\,(x,sign(y).|y|^{\frac{1}{\lambda_1}},
                sign(z).|z|^{\frac{1}{\lambda_2}},t,s).$
\vskip0.5cm
        $h_{4,3(\lambda)}(x,y,z,t,s)$: =\,

            \hskip2cm $=\,(x,sign(y).|y|^{\frac{1}{\lambda}},
            sign(z).|z|^{\frac{1}{\lambda}},t,s).$
\vskip0.5cm
        $h_{4,4(\lambda)}(x,y,z,t,s)$: =\,

                \hskip2cm $=\,(x,sign(y).|y|^{\frac{1}{\lambda}},z,t,s).$
\vskip0.5cm
        $h_{4,6(\lambda_1,\lambda_2)}(x,y,z,t,s)$: =
             $$=\, \left\{%
        \begin{array}{ll}
            (x, sign(y).|y|^{\frac{1}{\lambda_1}},
            sign(z).|z|^{\frac{1}{\lambda_2}},t,s-t.ln|t|), &\hbox{t $\neq$ 0;} \\
            (x, sign(y).|y|^{\frac{1}{\lambda_1}},
            sign(z).|z|^{\frac{1}{\lambda_2}},0,s), &\hbox{t = 0.} \\
        \end{array}%
        \right.$$
\vskip0.5cm
        $h_{4,7(\lambda)}(x,y,z,t,s)$: =
            $$=\, \left\{%
        \begin{array}{ll}
           (x, sign(y).|y|^{\frac{1}{\lambda}},
           sign(z).|z|^{\frac{1}{\lambda}},t,s-t.ln|t|), &\hbox{t $\neq$ 0;} \\
           (x, sign(y).|y|^{\frac{1}{\lambda}},
           sign(z).|z|^{\frac{1}{\lambda}},0,s), &\hbox{t = 0.} \\
        \end{array}%
        \right.$$

        $h_{4,8(\lambda)}(x,y,z,t,s)$: =
        $$=\, \left\{%
        \begin{array}{ll}
            (x, sign(y).|y|^{\frac{1}{\lambda}},sign(z-\frac{1}{\lambda}y.ln|y|\cdot
            \\\quad \cdot{|z-\frac{1}{\lambda}y.ln|y||}^{\frac{1}{\lambda}},t,s-t.ln|t|);
            &\hbox{y $\neq$ 0, t $\neq$ 0;} \\
            (x, 0,sign(z).|z|^{\frac{1}{\lambda}},t,s-t.ln|t|);
            &\hbox{y = 0, t $\neq$ 0;} \\
            (x, sign(y).|y|^{\frac{1}{\lambda}},sign(z-\frac{1}{\lambda}y.ln|y|\cdot
            \\\quad \cdot{|z-\frac{1}{\lambda}y.ln|y||}^{\frac{1}{\lambda}},0,s);
            &\hbox{y $\neq$ 0, t = 0;} \\
            (x, 0,sign(z).|z|^{\frac{1}{\lambda}},0,s); &\hbox{y = 0, t = 0.} \\
        \end{array}%
        \right. $$
\vskip0.5cm
        $h_{4,9(\lambda)}(x,y,z,t,s):=\,
        (\widetilde{x},\widetilde{y},\widetilde{z},\widetilde{t},\widetilde{s})$,
        \,where
\vskip 0.5cm
        $\begin{array}{l}
          \widetilde{x} = x ; \\
          \widetilde{y} = sign(y).|y|^{\frac{1}{\lambda}}; \\
          \widetilde{z} = z ; \\
          \widetilde{t} =\left\{%
            \begin{array}{ll}
              t-z.ln|z|, & \hbox{z $\neq0$;} \\
               t, & \hbox{z = 0;} \\
            \end{array}%
            \right.       \\
          \widetilde{s} = \left\{%
            \begin{array}{ll}
            s-\frac{1}{2}t.ln|z|-\\
            \quad -\frac{1}{2}(t-z.ln|z|).ln|t-z.ln|z||, & \hbox{$z \neq 0$, $t\neq{zln|z|}$;} \\
            s-\frac{1}{2}t.ln|z|, & \hbox{$z \neq 0, t = zln|z|$;} \\
            s, & \hbox{$z = 0, t = zln|z|$.} \\
            \end{array}%
            \right.      \\
        \end{array}$
\vskip1cm
        $h_{4,10(\lambda)}(x,y,z,t,s): =\,
        (\widetilde{x},\widetilde{y},\widetilde{z},\widetilde{t},\widetilde{s})$,
        \, where
\vskip0.5cm
        $\begin{array}{l}
         \widetilde{x} = x ; \\
         \widetilde{y} = y ; \\
         \widetilde{z} =\left\{%
            \begin{array}{ll}
            z-y.ln|y|, & \hbox{$y\neq0$;} \\
            z, & \hbox{$y = 0$;} \\
            \end{array}%
            \right. \\
         \widetilde{t} = \left\{%
            \begin{array}{ll}
            t-\frac{1}{2}z.ln|y| -\\
            \quad \frac{1}{2}(z-y.ln|y|).ln|z-y.ln|y||, & \hbox{$y \neq 0$, $z\neq{yln|y|}$;} \\
            t-\frac{1}{2}z.ln|y|, & \hbox{$y \neq 0, z = yln|y|$;} \\
            t, & \hbox{$y = 0, z = yln|y|$;} \\
            \end{array}%
            \right.     \\
           \widetilde{s} = \left\{%
            \begin{array}{ll}
            s-\frac{1}{3}t.ln|y|- \\
            \,\,\frac{1}{3}(t-\frac{1}{2}z.ln|y|).ln|y|+\\
            \,\,\frac{1}{3}(z.ln|y|-t-\frac{1}{2}y.ln^2|y|).ln|y|, & \hbox{$y \neq 0$;} \\
            s-\frac{1}{2}t.ln|z|- \\
            \,\,\frac{1}{2}(t-z.ln|z|).ln|t-z.ln|z||, & \hbox{$y=0\neq z$, $t\neq{zln|z|}$;} \\
            s-\frac{1}{2}t.ln|z|, & \hbox{$y=0\neq z, t = zln|z|$;} \\
            s, & \hbox{$y=0, z = 0$.} \\
            \end{array}%
            \right.     \\
        \end{array}$
\vskip1.5cm It is easy to verify that considered maps are
homeomorphisms which take leaves of each foliation listed in Set 1.1
of Theorem 3.2 except $(V,\,{\mathcal{F}}_{4,5})$ onto leaves of
this one. So these foliations are topologically equivalent to each
other.
        \item[1.2.] The topological equivalence of foliations in Set 1.2 of
        Theorem 3.2 is also verified similarly by considering homeomorphisms from
        $V\cong \mathbb{R}\times{(\mathbb{C}\times \mathbb{R}^2)}^{*}$ to
        oneself as follows

$h_{4,11(\lambda_1,\lambda_2,\varphi)} (x,re^{i\theta},t,s)$: =

\hskip1cm $ =\, (x,e^{(lnr+i\theta)(-ie^{i\varphi})},
sign(t).|t|^{\frac{1}{\lambda_1}},sign(s).|s|^{\frac{1}{\lambda_2}})$;
\vskip0.5cm
 $h_{4,12(\lambda,\varphi)} (x,re^{i\theta},t,s)$: =

\hskip1cm $=\,(x,e^{(lnr+i\theta)(-ie^{i\varphi})},
sign(t).|t|^{\frac{1}{\lambda}},sign(s).|s|^{\frac{1}{\lambda}})$;
\vskip0.5cm
$h_{4,13(\lambda,\varphi)} (x,re^{i\theta},t,s)$: =

\hskip1cm $ =\, \left\{%
\begin{array}{ll}
    (x, e^{(lnr+i\theta)(-ie^{i\varphi})}, sign(t).|t|^{\frac{1}{\lambda}},
    s - \frac{1}{\lambda}t.ln|t|), & \hbox{$t \neq 0$;} \\
    (x, e^{(lnr+i\theta)(-ie^{i\varphi})}, 0, s), & \hbox{$t = 0$.} \\
\end{array}%
\right. $
        \item[1.3.] Similarly, the homeomorphisms
        $$h_{4,14(\lambda,\mu,\varphi)}:V\cong\mathbb{R}
        \times{(\mathbb{C} \times \mathbb{C})}^{*}\rightarrow V$$

    $$(x,re^{i\theta},r'e^{i\theta'})\mapsto\left\{%
\begin{array}{ll}
    \Bigl(x, e^{(lnr+i\theta)(-ie^{i\varphi})}, \\
    \qquad e^{(lnr'+i\theta')(\frac{\mu}{\lambda^2 + \mu^2}
    -i \frac{\lambda}{\lambda^2 + \mu^2})}\Bigl ),
    & \hbox{$\lambda \neq 0$;} \\
    \Bigl(x, e^{(lnr+i\theta)(-ie^{i\varphi})}, \\
    \qquad e^{(lnr'+i\theta')\frac{1}{\mu}}\Bigl), & \hbox{$\lambda=0$.} \\
\end{array}%
\right.$$ take leaves of
$(V,\,{\mathcal{F}}_{4,14(\lambda,\mu,\varphi)})$ onto leaves of
$(V,\,{\mathcal{F}}_{4,14(0,1,\frac{\pi}{2})})$. Thus, foliations
listed in Set 1.3 of Theorem 3.2 are topologically equivalent to
each other.
    \end{enumerate}

    \item[2.] Now we prove Assertion 2 in
    Theorem 3.2. Firstly, it is easily seen that the following submersion:
    $$p_{4,5}:V\cong{(\mathbb{R}^4)}^*\times \mathbb{R}\cong\mathbf{S^3}
    \times \mathbb{R}_+\times \mathbb{R}\rightarrow \mathbf{S}^3;\,\,\,p_{4,5}(s,u,v): =\, s $$
    define Foliation $(V_4, {\mathcal{F}}_{4,5})$. Hence the
    foliations of ${\mathcal{F}}_1$ are trivial fibrations. Furthermore,
    the fibres of them are simply connected on $\mathbf{S}^3$.

    Nextly, let us consider the actions of the real commutative Lie group
    $\mathbb{R}^2$ on $V\cong{(\mathbb{R}^4)}^*\times \mathbb{R}$ as
    follows
    $$\rho_{4,12}: \mathbb{R}^2\times V \rightarrow V $$
    $$\rho_{4,12}((r,a),(x,y+iz,t,s))=(x+r,(y+iz).e^{-ia},te^a,se^a),$$
    where $(r,a)\in \mathbb{R}^2,(x,y+iz,t,s)\in V \cong \mathbb{R}
    \times{(\mathbb{C}\times \mathbb{R}^2)}^*$.
    $$\rho_{4,14}: \mathbb{R}^2\times V \rightarrow V  $$
    $$\rho_{4,14}((r,a),(x,y+iz,t+is))=(x+r,(y+iz).e^{-ia},(t+is)e^{-ia}),$$
    where $(r,a)\in \mathbb{R}^2,(x,y+iz,t+is)\in{V_4}\cong \mathbb{R}
    \times(\mathbb{C}\times \mathbb{C})^*$.

     It can be verified that the above actions
    $\rho_{4,12},\rho_{4,14}$ generate the foliations
    $(V,\,{\mathcal{F}}_{4,12(1,\frac{\pi}{2})}),\,
    (V,\,{\mathcal{F}}_{4,14(0,1,\frac{\pi}{2})})$ respectively. Hence,
    the foliations of ${\mathcal{F}}_2$, ${\mathcal{F}}_3$ can be
    given by suitable actions of $\mathbb{R}^2$ on the foliated
    manifolds $V\cong (\mathbb{R}^4)^* \times \mathbb{R}$.
    The proof is complete. \hfill{$\square$}
\end{description}

\subsection*{Concluding Remark}
We close the paper with some remarks as follows.
\begin{itemize}
    \item The topological structure of considered MD5-foliations
    will be characterized profoundly when we study the Connes'
    C*-algebras of them. In the next paper we shall be concerned
    with this problem.
    \item It should be noted that if the foliation $(V,\,\mathcal{F})$
    comes from a fibration  $p: V\rightarrow B$ (with connected fibers)
    then the Connes' C*-algebra $C^*(V,\,\mathcal{F})$ is isomorphic to the
    tensor product $C_0(B) \otimes \mathcal{K}$, where $C_0(B)$ is
    the algebra of $\mathbb{C}$-values continuous functions on B vanishing
    at infinity and $\mathcal{K}$ denotes the C*-algebra of compact
    operators on an (infinite dimensional) separable Hilbert space
    (see [1, Section 5]). So we have the following assertion as one
    direct consequence of Theorem 3.2.
\end{itemize}
\vskip0.5cm
\begin{corollary}(The C*- algebras of MD5-foliations
of the type ${\mathcal{F}}_1$)

The Connes's C*- algebras of all MD5-foliations of the type
${\mathcal{F}}_1$ are isomorphic to the C*-algebra
$C(\mathbf{S}^3)\otimes \mathcal{K}$. \hfill{$\square$}

\end{corollary}
\subsection*{Acknowledgement}
 The authors would like take this opportunity to thank
Prof. DSc. Do Ngoc Diep for his excellent advice and support. They
wish to thank Prof. Nguyen Van Sanh for his encouragement. Thanks
are due also to the Scientific Committee and Organizing Committee of
The Second International Congress In Algebras and Combinatorics -
July 2007, Xi'an, China for inviting the first author to come and
give talk on this topic at the congress.

\end{document}